\pdfoutput=1

\documentclass[11pt,reqno,final]{amsart}
\usepackage{amssymb,amsmath,color,tikz-cd,stmaryrd}
\usepackage{geometry}
\newtheorem{theorem}{Theorem}[section]
\newtheorem{lemma}[theorem]{Lemma}
\newtheorem{proposition}[theorem]{Proposition}

\newtheorem{example}{\it Example\/}
\newtheorem*{formality}{Formality theorem}
\newtheorem*{kontsevichduflo}{Kontsevich--Duflo type theorem}

\newcommand{\CE}{\operatorname{CE}}
\newcommand{\hkr}{\operatorname{hkr}}
\newcommand{\tot}{\operatorname{tot}}
\newcommand{\poly}{\operatorname{poly}}
\newcommand{\End}{\operatorname{End}}
\newcommand{\Todd}{\operatorname{Td}}
\newcommand{\cTodd}{\operatorname{td}}
\newcommand{\Ext}{\operatorname{Ext}}
\newcommand{\id}{\operatorname{id}}

\newcommand{\gTpoly}[1]{\tot\big(\Lambda^{#1} \mathfrak{g}\dual \otimes_{\KK} \Tpoly{#1} \big)}
\newcommand{\gDpoly}[1]{\tot\big(\Lambda^{#1} \mathfrak{g}\dual \otimes_{\KK} \Dpoly{#1} \big)}
\newcommand{\Tpoly}[1]{T_{\poly}^{#1}(M)}
\newcommand{\Dpoly}[1]{D_{\poly}^{#1}(M)}
\newcommand{\sections}[1]{\Gamma(#1)}
\newcommand{\schouten}[2]{[#1,#2]} 
\newcommand{\gerstenhaber}[2]{\llbracket #1,#2\rrbracket} 
\newcommand{\hochschild}{d_{H}}
\newcommand{\into}{\hookrightarrow}
\newcommand{\iso}{\xto{\cong}}
\newcommand{\dual}{^{\vee}}
\newcommand{\KK}{\Bbbk} 
\newcommand{\atiyahcocycle}{{R^\nabla_{1,1}}}
\newcommand{\atiyahclass}{\alpha}
\newcommand{\cO}{\mathcal{O}}
\newcommand{\xto}[1]{\xrightarrow{#1}}
\newcommand{\RR}{\mathbb{R}} 
\newcommand{\CC}{\mathbb{C}} 
\newcommand{\vecf}[1]{\mathfrak{X}(#1)}
\newcommand{\Tpolycmplx}{\Tpoly{\bullet}\xto{0}\Tpoly{\bullet+1}}
\newcommand{\Dpolycmplx}{\Dpoly{\bullet}\xto{\hochschild}\Dpoly{\bullet+1}}
\newcommand{\liederivative}[1]{\mathcal{L}_{#1}}

\begin{document}

\title{Formality Theorem for $\mathfrak{g}$-manifolds}

\dedicatory{\`A la m\'emoire de notre ami Jacky Mich\'ea}

\thanks{Research partially supported by NSF grants DMS-1406668 and DMS-1101827, and NSA grant H98230-14-1-0153.}

\author{Hsuan-Yi Liao}
\address{Department of Mathematics, Pennsylvania State University}
\email{hul170@psu.edu}

\author{Mathieu Sti\'enon}
\address{Department of Mathematics, Pennsylvania State University}
\email{stienon@psu.edu} 

\author{Ping Xu}
\address{Department of Mathematics, Pennsylvania State University}
\email{ping@math.psu.edu}

\begin{abstract}
To any $\mathfrak{g}$-manifold $M$ are associated two dglas
$\gTpoly\bullet$ and $\gDpoly\bullet$, whose cohomologies 
$H^\bullet_{\CE}\big(\mathfrak{g},\Tpolycmplx\big)$
and $H^\bullet_{\CE}\big(\mathfrak{g},\Dpolycmplx\big)$ are Gerstenhaber algebras. 
We establish a formality theorem for $\mathfrak{g}$-manifolds: 
there exists an $L_\infty$ quasi-isomorphism 
$\Phi: \gTpoly\bullet\to \gDpoly\bullet$
whose first `Taylor coefficient' (1) is equal to the 
Hochschild--Kostant--Rosenberg map twisted by the square root of the Todd 
cocycle of the $\mathfrak{g}$-manifold $M$ and (2) induces an isomorphism of 
Gerstenhaber algebras on the level of cohomology.
Consequently, the Hochschild--Kostant--Rosenberg map twisted by 
the square root of the Todd class of the $\mathfrak{g}$-manifold $M$ 
is an isomorphism of Gerstenhaber algebras from 
$H^\bullet_{\CE}\big(\mathfrak{g},\Tpolycmplx\big)$ 
to $H^\bullet_{\CE}\big(\mathfrak{g},\Dpolycmplx\big)$.
\end{abstract}

\maketitle

\section{Introduction}

Two differential graded Lie algebras (dglas) are canonically associated 
with a given smooth manifold $M$: the dgla of polyvector fields 
$\Tpoly\bullet=\bigoplus_{k=-1}^\infty \sections{\Lambda^{k+1}T_M}$, which
is endowed with the zero differential and the Schouten bracket 
$\schouten{\;}{\;}$, and the dgla of polydifferential operators 
$\Dpoly\bullet= \bigoplus_{k=-1}^\infty \Dpoly{k}$, which is endowed 
with the Hochschild differential $\hochschild$ and the Gerstenhaber bracket 
$\gerstenhaber{\;}{\;}$.
Here $\Dpoly{-1}$ denotes the algebra of smooth functions $\mathcal{R}=C^\infty(M)$,
$\Dpoly{0}$ the algebra of differential operators on $M$, 
and $\Dpoly{k}$ (with $k\geq 0$) the space of $(k+1)$-differential operators 
on $M$, i.e.\ the tensor product 
$\Dpoly{0}\otimes_{\mathcal{R}}\cdots\otimes_{\mathcal{R}}\Dpoly{0}$ 
of $(k+1)$ copies of the left $\mathcal{R}$-module $\Dpoly{0}$.
The classical Hochschild--Kostant--Rosenberg (HKR) theorem 
\cite{MR0142598,MR2062626} states that 
the Hochschild--Kostant--Rosenberg map, 
the natural embedding $\hkr:\Tpoly\bullet \into \Dpoly\bullet$ 
defined by Equation~\eqref{eq:hkr}, determines 
an isomorphism of Gerstenhaber algebras 
$\hkr: \Tpoly\bullet \iso H^\bullet(\Dpoly\bullet, \hochschild)$ 
on the cohomology level --- the products on $\Tpoly\bullet$ and 
$\Dpoly\bullet$ are the wedge product and the cup product respectively. 
However, the HKR map $\hkr:\Tpoly\bullet\into\Dpoly\bullet$
is not a morphism of dglas. Kontsevich's celebrated
formality theorem states that the HKR map $\hkr$ extends
to an $L_\infty$ quasi-isomorphism from $\Tpoly\bullet$
to $\Dpoly\bullet$ \cite{MR2062626,MR2699544}. 
The formality theorem is highly non trivial and has many applications, 
one of which is the deformation quantization of Poisson manifolds. 

In this Note, we study the Gerstenhaber algebra structures associated with 
a $\mathfrak{g}$-manifold and we establish a formality theorem for 
$\mathfrak{g}$-manifolds.
By a $\mathfrak{g}$-manifold, we mean a  smooth manifold equipped with 
an infinitesimal action of a Lie algebra $\mathfrak{g}$.
In this situation, the analogues of $\Tpoly\bullet$ and $\Dpoly\bullet$
are the Chevalley--Eilenberg complexes 
$\tot\big(\Lambda^\bullet\mathfrak{g}\dual\otimes_{\KK}\Tpoly\bullet\big)$
and $\tot\big(\Lambda^\bullet\mathfrak{g}\dual\otimes_{\KK}\Dpoly\bullet\big)$, respectively 
--- they are briefly mentioned in Dolgushev's work~\cite[concluding remarks]{MR2102846}. 
Both of them are naturally dglas (see Lemma~\ref{lem:3.1} and Lemma~\ref{lem:3.2}) and their cohomologies are Gerstenhaber algebras.

In order to state the formality theorem and the precise relation between 
these two Gerstenhaber algebras, one must take into consideration 
the obstruction to the existence of a $\mathfrak{g}$-invariant affine connection 
on $M$, the Atiyah cocycle 
$\atiyahcocycle\in\mathfrak{g}\dual\otimes\sections{T_M\dual\otimes\End T_M}$,
which is a Chevalley--Eilenberg $1$-cocycle of the $\mathfrak{g}$-module
$\sections{T_M\dual \otimes \End T_M}$. 
More precisely, we must call upon its cohomology class, the Atiyah  class $\atiyahclass_{M/\mathfrak{g}}\in H_{\CE}^1(\mathfrak{g}, 
\sections{T_M\dual\otimes\End T_M})$, which we introduce in Proposition~\ref{pro:Atiyah}. 

The Todd cocycle $\cTodd_{M/\mathfrak{g}}\in\bigoplus_{k=0} \Lambda^{k}\mathfrak{g}^\vee \otimes\Omega^k(M)$ 
of a $\mathfrak{g}$-manifold $M$ is defined in terms of the Atiyah cocycle 
in Equation~\eqref{eq:todd}. 
The corresponding class in Chevalley--Eilenberg cohomology 
is the Todd class 
$\Todd_{M/\mathfrak{g}}\in\bigoplus_{k=0}H_{\CE}^{k}(\mathfrak{g},\Omega^k (M))$.
See Equation~\eqref{eq:Todd}.

The main results of this Note are a formality theorem for $\mathfrak{g}$-manifolds 
and its consequence: a Kontsevich--Duflo type theorem for $\mathfrak{g}$-manifolds. 

\begin{formality}
Given a $\mathfrak{g}$-manifold $M$ 
and an affine torsionfree connection $\nabla$ on $M$,
there exists an $L_\infty$ quasi-isomorphism $\Phi$ 
from the dgla
$\tot\big(\Lambda^\bullet \mathfrak{g}\dual \otimes_{\KK} \Tpoly\bullet\big)$
to the dgla
$\tot \big(\Lambda^\bullet\mathfrak{g}\dual \otimes_{\KK} \Dpoly\bullet\big)$
whose first `Taylor coefficient' $\Phi_1$ satisfies the following two properties: 
\begin{enumerate} 
\item $\Phi_1$ is, up to homotopy, 
an isomorphism of associative algebras (and hence induces an isomorphism of associative algebras of the homologies); 
\item $\Phi_1$ is equal to the composition 
$\hkr\circ\cTodd^{\frac{1}{2}}_{M/\mathfrak{g}}$ 
of the HKR map and the action of the square root of the Todd cocycle 
$\cTodd^{\frac{1}{2}}_{M/\mathfrak{g}} \in \bigoplus_{k=0} \Lambda^{k}\mathfrak{g}^\vee \otimes \Omega^k (M)$ on
$\tot\big(\Lambda^\bullet\mathfrak{g}\dual\otimes_{\KK}\Tpoly\bullet\big)$
by contraction. 
\end{enumerate}
\end{formality}

\begin{kontsevichduflo}
Given a $\mathfrak{g}$-manifold $M$, the map
\[ \hkr\circ\Todd^{\frac{1}{2}}_{M/\mathfrak{g}}: 
H^\bullet_{\CE}\big(\mathfrak{g},\Tpolycmplx\big)
\longrightarrow H^\bullet_{\CE}\big(\mathfrak{g},\Dpolycmplx\big) \]
is an isomorphism of Gerstenhaber algebras.
Here $H_{\CE}^k(\mathfrak{g},E^\bullet\xto{d_E}E^{\bullet+1})$ denotes the 
Chevalley--Eilenberg cohomology of $\mathfrak{g}$ with coefficients in the complex of $\mathfrak{g}$-modules $E^\bullet$. 
It is understood that the square root $\Todd^{\frac{1}{2}}_{M/\mathfrak{g}}$ of the 
Todd class $\Todd_{M/\mathfrak{g}}\in\bigoplus_{k=0} 
H_{\CE}^{k}\big(\mathfrak{g}, \Omega^k (M)\big)$ acts on 
$H^\bullet_{\CE}\big(\mathfrak{g},\Tpolycmplx\big)$ by contraction. 
\end{kontsevichduflo}

The theorem above is parallel in spirit to an analogue of Duflo's Theorem --- a classical result of Lie theory --- 
discovered by Kontsevich in complex geometry \cite{MR2062626}. Kontsevich observed that, for a complex manifold $X$,
the composition $\hkr\circ(\Todd_X)^{\frac{1}{2}}:
H^\bullet(X,\Lambda^\bullet T_X)\xto{\cong} HH^\bullet(X)$ 
is an isomorphism of associative algebras. 
Here $\Todd_X$ denotes the Todd class of the tangent bundle $T_X$ 
and $HH^\bullet(X)$ denotes the Hochschild cohomology groups 
of the complex manifold $X$, i.e.\ the groups 
$\Ext_{\cO_{X\times X}}^\bullet(\cO_{\Delta},\cO_{\Delta})$.
The multiplications on $H^\bullet(X,\Lambda^{\bullet}T_X)$
and $HH^{\bullet}(X)$ are the wedge product and the Yoneda product respectively.
A detailed proof of Kontsevich's result appeared in~\cite{MR2646112}. 
It is worth mentioning that the map $\hkr\circ(\Todd_X)^{\frac{1}{2}}$ 
actually respects the Gerstenhaber algebra structures 
on both cohomologies; this was brought to light in~\cite{MR2646112}.

\section{Preliminary: Chevalley--Eilenberg cohomology}

Let $\mathfrak{g}$ be a Lie algebra over $\KK$ ($\KK$ is $\RR$ or $\CC$). 
Given a $\mathfrak{g}$-module $E$, one may consider 
the Chevalley--Eilenberg cochain complex 
\[ \begin{tikzcd}
\cdots \ar{r}& \Lambda^{p-1}\mathfrak{g}\dual \otimes_{\KK} E \ar{r}{d_{\CE}} & \Lambda^{p}\mathfrak{g}\dual \otimes_{\KK} E \ar{r}{d_{\CE}} & \Lambda^{p+1}\mathfrak{g}\dual \otimes_{\KK} E \ar{r} & \cdots,
\end{tikzcd} \]
where $d_{\CE}$ is the Chevalley--Eilenberg differential.
More generally, given a bounded below complex of left $\mathfrak{g}$-modules
\[ \begin{tikzcd}
\cdots \ar{r}& E^{p-1} \ar{r}{d_E} & E^p \ar{r}{d_E} & E^{p+1} \ar{r} & \cdots 
,\end{tikzcd} \]
we may consider the double complex: 
\[ \begin{tikzcd}[ampersand replacement=\&]
\& \vdots \& \vdots \& \vdots \& \\ 
\cdots \arrow{r} \& 
\Lambda^{p-1}\mathfrak{g}\dual \otimes_{\KK} E^{q+1} 
\arrow{r}{d_{\CE}} \arrow{u} \& 
\Lambda^{p}\mathfrak{g}\dual\otimes_{\KK}E^{q+1} 
\arrow{r}{d_{\CE}} \arrow{u} \& 
\Lambda^{p+1}\mathfrak{g}\dual\otimes_{\KK}E^{q+1} 
\arrow{r} \arrow{u} \& \cdots \\ 
\cdots \arrow{r} \& 
\Lambda^{p-1}\mathfrak{g}\dual\otimes_{\KK}E^{q} 
\arrow{r}[swap]{d_{\CE}} \arrow{u}{(-1)^{p-1}\id\otimes  d_{E}} \& 
\Lambda^{p}\mathfrak{g}\dual\otimes_{\KK}E^{q} 
\arrow{r}[swap]{d_{\CE}} \arrow{u}{(-1)^{p}\id\otimes d_{E}} \& 
\Lambda^{p+1}\mathfrak{g}\dual\otimes_{\KK}E^{q} 
\arrow{r} \arrow{u}{(-1)^{p+1}\id\otimes d_{E}} \& \cdots \\ 
\& \vdots \arrow{u} \& \vdots \arrow{u} \& \vdots \arrow{u} \& 
\end{tikzcd} \]
where $d_{\CE}$ is the Chevalley--Eilenberg differential corresponding
to the $\mathfrak{g}$-module structure on $E^\bullet$.
By definition, the Chevalley--Eilenberg cohomology of $\mathfrak{g}$ with coefficients in 
the complex of $\mathfrak{g}$-modules $(E^\bullet, d_E)$ is the total cohomology 
of the double complex above:
\[ H_{\CE}^k(\mathfrak{g},E^\bullet\xto{d_E}E^{\bullet+1})=
H^k\big(\tot(\Lambda^\bullet\mathfrak{g}\dual \otimes_{\KK} E^\bullet)
\big)\]

\section{Hochschild--Kostant--Rosenberg theorem for $\mathfrak{g}$-manifolds}\label{HKR}

\subsection{Polyvector fields}

Let $M$ be a $\mathfrak{g}$-manifold with infinitesimal action given by
a Lie algebra morphism $\varphi: \mathfrak{g} \to \vecf{M}$.
It is well known that the space of polyvector fields 
 $\Tpoly\bullet=\bigoplus_{k=-1}^\infty \sections{\Lambda^{k+1}T_M}$
on $M$, together with the wedge product and the Schouten bracket 
$\schouten{\;}{\;}$, forms a Gerstenhaber algebra. 
Moreover, the $\mathfrak{g}$-action on $M$ and the Schouten bracket together determine a $\mathfrak{g}$-module structure on $\Tpoly{k}$ for each $k\geq -1$: 
\[ a \cdot \gamma = \schouten{\varphi(a)}{\gamma}  \qquad \forall \; a 
\in \mathfrak{g}, \; \gamma \in \Tpoly{k} .\]

Therefore $\cdots\to\Tpoly{k}\xto{0}\Tpoly{k+1}\to\cdots$ is a 
complex of $\mathfrak{g}$-modules. Its Chevalley--Eilenberg cohomology 
\[ H_{\CE}^k \big(\mathfrak{g}, \Tpoly{\bullet}\xto{0}\Tpoly{\bullet+1}\big) = H^k\big(\tot(\Lambda^\bullet\mathfrak{g}\dual \otimes_{\KK} \Tpoly\bullet )\big) \]
is the total cohomology of the double complex:
\[ \begin{tikzcd}[ampersand replacement=\&]
\& \vdots \& \vdots \& \vdots \& \\ 
\cdots \arrow{r} \& 
\Lambda^{p-1}\mathfrak{g}\dual \otimes_{\KK} \Tpoly{q+1} 
\arrow{r}{d_{\CE}} \arrow{u}{0} \& 
\Lambda^{p}\mathfrak{g}\dual\otimes_{\KK} \Tpoly{q+1}
\arrow{r}{d_{\CE}} \arrow{u}{0} \& 
\Lambda^{p+1}\mathfrak{g}\dual\otimes_{\KK} \Tpoly{q+1} 
\arrow{r} \arrow{u}{0} \& \cdots \\ 
\cdots \arrow{r} \& 
\Lambda^{p-1}\mathfrak{g}\dual\otimes_{\KK}\Tpoly{q}
\arrow{r}[swap]{d_{\CE}} \arrow{u}{0} \& 
\Lambda^{p}\mathfrak{g}\dual\otimes_{\KK}\Tpoly{q}
\arrow{r}[swap]{d_{\CE}} \arrow{u}{0} \& 
\Lambda^{p+1}\mathfrak{g}\dual\otimes_{\KK}\Tpoly{q} 
\arrow{r} \arrow{u}{0} \& \cdots \\ 
\& \vdots \arrow{u}{0} \& \vdots \arrow{u}{0} \& \vdots \arrow{u}{0} \& 
\end{tikzcd} \]

Extend the Schouten bracket $\schouten{\;}{\;}$ on $\Tpoly\bullet$
to $\Lambda^\bullet \mathfrak{g}\dual \otimes_{\KK} \Tpoly\bullet$ as follows:
\begin{equation}\label{bracketTpoly}
\schouten{\alpha \otimes \mathcal{X}}{\beta \otimes \mathcal{Y}} = (-1)^{q_1 p_2} \alpha \wedge \beta \otimes \schouten{\mathcal{X}}{\mathcal{Y}}
\end{equation}
for any $\alpha \otimes \mathcal{X} \in \Lambda^{p_1}\mathfrak{g}\dual \otimes_{\KK} \Tpoly{q_1}$ and $\beta \otimes \mathcal{Y} \in \Lambda^{p_2} \mathfrak{g}\dual \otimes_{\KK} \Tpoly{q_2}$.

The following lemma can be easily verified.
\begin{lemma}
\label{lem:3.1}
The graded $\KK$-vector space 
$\tot\big(\Lambda^\bullet \mathfrak{g}\dual \otimes_{\KK} \Tpoly\bullet\big)$,
together with the Chevalley-Eilenberg differential $d_{\CE}$,
the wedge product $\wedge$ and the bracket defined by Equation~\eqref{bracketTpoly} is
a differential Gerstenhaber algebra. As a consequence, 
$H^\bullet_{\CE}(\mathfrak{g},\Tpolycmplx)$ is a Gerstenhaber algebra.
\end{lemma}

\subsection{Polydifferential operators}

On a smooth  manifold $M$, one also has the dgla 
of polydifferential operators, $\Dpoly\bullet$. 

Let $M$ be a manifold, let $\mathcal{R}$ denote its algebra of smooth functions $C^\infty(M)$, 
and let $\Dpoly{0}$ denote the algebra of differential operators on $M$. 
Denote by $\Dpoly{k}$, $k\geq 0$, the space of $(k+1)$-differential
 operators on $M$, i.e.\ the tensor product $\Dpoly{0} \otimes_{\mathcal{R}} 
\cdots\otimes_{\mathcal{R}}\Dpoly{0}$ of $(k+1)$ copies of the left $\mathcal{R}$-module
 $\Dpoly{0}$. Denote also by $\Dpoly{-1}$ the space of smooth functions $\mathcal{R}=C^\infty(M)$.
It is well known that endowing $\Dpoly\bullet= \bigoplus_{k=-1}^\infty \Dpoly{k}$ 
with the Hochschild differential $\hochschild$,
the cup product $\Dpoly{k}\otimes\Dpoly{l}\xto{\smile}\Dpoly{k+l+1}$, 
and the Gerstenhaber bracket $\gerstenhaber{\;}{\;}$ makes it a Gerstenhaber 
algebra up to homotopy \cite{MR0161898}.

Following our earlier notations, now assume that $M$ is a 
$\mathfrak{g}$-manifold with infinitesimal action 
$\varphi: \mathfrak{g} \to \vecf{M}$.
Analogously to the polyvector field case, 
the Lie algebra $\mathfrak{g}$ acts on $\Dpoly\bullet$ by:
\[ a \cdot \mu =\gerstenhaber{\varphi(a)}{\mu}  \qquad \forall \; a \in \mathfrak{g}, \; \mu \in \Dpoly\bullet .\] 
Since the Gerstenhaber bracket satisfies the graded Jacobi identity, 
this infinitesimal $\mathfrak{g}$-action on $\Dpoly{\bullet}$ 
is compatible with the Hochschild differential. 
Consequently $\cdots \to \Dpoly{k}\xto{\hochschild} \Dpoly{k+1}\to\cdots$
is a complex of $\mathfrak{g}$-modules, and therefore we have the 
Chevalley--Eilenberg cohomology
\[ H_{\CE}^k \big(\mathfrak{g},\Dpolycmplx\big) = H^k\big(\tot(\Lambda^\bullet\mathfrak{g}\dual \otimes_{\KK} \Dpoly\bullet)\big) ,\]
which is, by definition, the total cohomology of the double complex 
\[ \begin{tikzcd}[ampersand replacement=\&]
\& \vdots \& \vdots \& \vdots \& \\ 
\cdots \arrow{r} \& 
\Lambda^{p-1}\mathfrak{g}\dual \otimes_{\KK} \Dpoly{q+1} 
\arrow{r}{d_{\CE}} \arrow{u}{(-1)^{p-1}\id\otimes\hochschild} \& 
\Lambda^{p}\mathfrak{g}\dual\otimes_{\KK}\Dpoly{q+1}
\arrow{r}{d_{\CE}} \arrow{u}{(-1)^{p}\id\otimes\hochschild} \& 
\Lambda^{p+1}\mathfrak{g}\dual\otimes_{\KK} \Dpoly{q+1} 
\arrow{r} \arrow{u}{(-1)^{p+1}\id\otimes\hochschild} \& \cdots \\ 
\cdots \arrow{r} \& 
\Lambda^{p-1}\mathfrak{g}\dual\otimes_{\KK}\Dpoly{q}
\arrow{r}[swap]{d_{\CE}} \arrow{u}{(-1)^{p-1}\id\otimes\hochschild} \& 
\Lambda^{p}\mathfrak{g}\dual\otimes_{\KK}\Dpoly{q}
\arrow{r}[swap]{d_{\CE}} \arrow{u}{(-1)^{p}\id\otimes\hochschild} \& 
\Lambda^{p+1}\mathfrak{g}\dual\otimes_{\KK}\Dpoly{q} 
\arrow{r} \arrow{u}{(-1)^{p+1}\id\otimes\hochschild} \& \cdots \\ 
\& \vdots \arrow{u}{(-1)^{p-1}\id\otimes\hochschild} \& 
\vdots \arrow{u}{(-1)^{p}\id\otimes\hochschild} 
\& \vdots \arrow{u}{(-1)^{p+1}\id\otimes\hochschild} \& 
\end{tikzcd} \]

Extend the cup product $\smile$ and the Gerstenhaber bracket $\gerstenhaber{\;}{\;}$
to $\Lambda^\bullet \mathfrak{g}\dual \otimes_{\KK} \Dpoly\bullet$ as follows:
\begin{align}
\label{cup product}
(\alpha \otimes \xi) \smile (\beta \otimes \eta) & = (-1)^{q_1 p_2} (\alpha \wedge \beta) \otimes (\xi \smile \eta)\\
\label{bracketDpoly}
\gerstenhaber{\alpha \otimes \xi}{ \beta \otimes \eta} & = (-1)^{q_1 p_2} \alpha \wedge \beta \otimes \gerstenhaber{\xi}{\eta} 
\end{align}
for any $\alpha \otimes \xi \in \Lambda^{p_1}\mathfrak{g}\dual \otimes_{\KK} \Dpoly{q_1}$ and $\beta \otimes \eta \in \Lambda^{p_2} \mathfrak{g}\dual \otimes_{\KK} \Dpoly{q_2}$.

Again the following lemma is immediate.

\begin{lemma}
\label{lem:3.2}
\begin{enumerate}
\item The graded $\KK$-vector space $\tot(\Lambda^\bullet\mathfrak{g}\dual \otimes_{\KK} \Dpoly\bullet$, together with the differential
$d_{\CE}+\id\otimes\hochschild$ and the Gerstenhaber bracket 
$\gerstenhaber{\cdot}{\cdot}$ defined by Equation~\eqref{bracketDpoly},
is a dgla.
\item The graded $\KK$-vector space 
$H^\bullet_{\CE}(\mathfrak{g},\Dpolycmplx)$,
together with the cup product and the Gerstenhaber bracket 
defined by Equations~\eqref{cup product} and~\eqref{bracketDpoly}, 
is a Gerstenhaber algebra.
\end{enumerate}
\end{lemma}

\subsection{Hochschild--Kostant--Rosenberg theorem}

Given a smooth manifold $M$, there is a natural embedding $\hkr:\Tpoly\bullet \into \Dpoly\bullet$, called Hochschild--Kostant--Rosenberg map, and defined by 
\begin{equation}
\label{eq:hkr}
 \hkr(X_1 \wedge \cdots \wedge X_k)= \frac{1}{k!} \sum_{\sigma \in S_k} \operatorname{sgn}(\sigma) X_{\sigma(1)} \otimes \cdots \otimes X_{\sigma(k)}, \ \ \forall X_i\in\vecf{M},
\end{equation}
where $S_k$ is the symmetric group on $k$ objects. 
The Hochschild--Kostant--Rosenberg theorem for smooth manifolds states that
$\hkr$ is a quasi-isomorphism, i.e.\ the induced morphism in cohomology 
$\hkr: \Tpoly\bullet \xto{\cong} H^\bullet(\Dpoly\bullet, \hochschild)$ 
is an isomorphism of vector spaces \cite{MR0142598,MR2062626}. 

Suppose we are given a $\mathfrak{g}$-manifold $M$. 
Then the map 
$\id \otimes \hkr: \Lambda^\bullet \mathfrak{g}\dual \otimes_{\KK} \Tpoly\bullet 
\to \Lambda^\bullet \mathfrak{g}\dual \otimes_{\KK} \Dpoly\bullet$ 
is a morphism of double complexes. 
Abusing notations, the induced morphism on Chevalley--Eilenberg cohomologies 
will also be denoted by $\hkr$. 

\begin{proposition}[\cite{LSX:16}]
Let $M$ be a $\mathfrak{g}$-manifold. 
The Hochschild--Kostant--Rosenberg map
\[ \hkr: H^\bullet_{\CE}(\mathfrak{g},\Tpolycmplx)
\xto{\cong} H^\bullet_{\CE}(\mathfrak{g},\Dpolycmplx) \] 
is an isomorphism of vector spaces.
\end{proposition}

The proof is a straightforward spectral sequence computation relying 
on the classical Hochschild--Kostant--Rosenberg theorem for smooth manifolds.

\section{Atiyah class of a $\mathfrak{g}$-manifold}

The Atiyah class was originally introduced by Atiyah for holomorphic vector bundles \cite{MR0086359}. Atiyah classes can also be defined for Lie algebroid pairs \cite{MR3439229} and dg vector bundles \cite{MR3319134}.
In this section, we introduce the notions of Atiyah class and Todd class of a $\mathfrak{g}$-manifold. 

Let $M$ be a $\mathfrak{g}$-manifold with infinitesimal action 
$\mathfrak{g}\ni a\mapsto\hat{a}\in\vecf{M}$. 
Given an affine connection $\nabla$ on $M$, the Atiyah 1-cocycle associated 
with $\nabla$ is defined as the map 
$\atiyahcocycle: \mathfrak{g} \times \vecf{M} \to \End_{\mathcal{R}} \vecf{M}$ given by
\[ \atiyahcocycle(a,X)=\liederivative{\hat{a}}\circ\nabla_X - \nabla_X\circ\liederivative{\hat{a}} 
- \nabla_{\liederivative{\hat{a}}X} ,\]
where $a\in\mathfrak{g}$, $X\in\vecf{M}$, and $\mathcal{R}=C^\infty(M)$.

Following \cite{MR3439229}, we  prove the following

\begin{proposition} 
\label{pro:Atiyah}
\begin{enumerate}
\item The Atiyah cocycle $\atiyahcocycle \in \mathfrak{g}\dual \otimes \sections{T_M\dual \otimes \End T_M}$ is a Chevalley--Eilenberg $1$-cocycle of the $\mathfrak{g}$-module $\sections{T_M\dual \otimes \End T_M}$.
\item The cohomology class $\atiyahclass_{M/\mathfrak{g}}\in 
H_{\CE}^1(\mathfrak{g}, \sections{T_M\dual \otimes \End T_M})$ of the $1$-cocycle $\atiyahcocycle$
does not depend on the choice of connection $\nabla$.
\end{enumerate}
\end{proposition}

The cohomology class $\atiyahclass_{M/\mathfrak{g}}$ is called 
 the Atiyah class of the $\mathfrak{g}$-manifold $M$. It is
the obstruction class to the existence of a $\mathfrak{g}$-invariant connection on $M$, i.e.\ an affine connection $\nabla$ on $M$  satisfying
\[ [\hat{a},\nabla_X Y] = \nabla_{[\hat{a},X]} Y + \nabla_X [\hat{a},Y] \] 
for all $a\in\mathfrak{g}$ and $X,Y\in\vecf{M}$. 

\begin{proposition}
Let $M$ be a $\mathfrak{g}$-manifold. 
The Atiyah class $\atiyahclass_{M/\mathfrak{g}}$ of $M$ vanishes 
if and only if there exists a $\mathfrak{g}$-invariant connection on $M$.
\end{proposition}

Note that if $\mathfrak{g}$ is a compact Lie algebra, 
$\atiyahclass_{M/\mathfrak{g}}$  
vanishes since $\mathfrak{g}$-invariant connections always exist.

The Todd class of  complex vector bundles  plays an important role 
in the Riemann--Roch theorem. 
In our context, the Todd cocycle of a $\mathfrak{g}$-manifold $M$ 
is the   Chevalley--Eilenberg cocycle
\begin{equation}
\label{eq:todd}
\cTodd_{M/\mathfrak{g}} =\det\left(\frac{\atiyahcocycle}
{1-e^{-\atiyahcocycle}} \right)
\in \bigoplus_{k=0} \Lambda^{k}\mathfrak{g}^\vee \otimes \Omega^k (M),
\end{equation}
with $\Omega^k (M)$, $k\geq 0$,  being the natural $\mathfrak{g}$-module.
Its corresponding Chevalley--Eilenberg cohomology class is the
{\em Todd class} $\Todd_{M/\mathfrak{g}}$.  Alternatively
\begin{equation}
\label{eq:Todd}
 \Todd_{M/\mathfrak{g}} =\det\left(\frac{\atiyahclass_{M/\mathfrak{g}}}
{1-e^{-\atiyahclass_{M/\mathfrak{g}}}}\right)
\in \bigoplus_{k=0} H_{\CE}^{k}(\mathfrak{g}, \Omega^k (M)) .
\end{equation}
Since the Lie algebra $\mathfrak{g}$ is finite dimensional, 
the above expression for the Todd class $\Todd_{M/\mathfrak{g}}$ 
reduces to a finite sum.

\begin{example}
Consider the case of the $1$-dimensional abelian Lie algebra $\mathfrak{g}= \RR$ 
acting on the real line $M=\RR$.
The infinitesimal action is uniquely determined by a
vector field $Q=q(x)\frac{d}{d x}\in\vecf{\RR}$. 
The Chevalley--Eilenberg complex 
$\big(\Lambda^\bullet\mathfrak{g}\dual\otimes 
\sections{T_M\dual\otimes\End T_M},d_{\CE}\big)$ 
is then isomorphic to the 2-term  complex
\[ \begin{tikzcd}
0 \ar{r} & C^\infty(\RR) \ar{r}{d_Q} & C^\infty(\RR) \ar{r} & 0,
\end{tikzcd} \] 
where the map $d_Q$ is given by
\[ d_Q(f) = \frac{d(fq)}{dx} = f'q+fq' ,\]
for $f\in C^\infty(\RR)$. 
Let $\nabla$ be the  trivial affine connection on the manifold $M=\RR$, 
i.e.\ $\nabla_{\frac{d}{dx}}\frac{d}{dx}=0$. Under the above isomorphism, 
the Atiyah 1-cocycle $\atiyahcocycle$ is simply the second order derivative of $q$: 
\[ \atiyahcocycle = q'' \in C^\infty(\RR) \cong 
\mathfrak{g}\dual \otimes \sections{T_M\dual \otimes \End T_M} .\] 
As a consequence, the Atiyah class vanishes if and only if 
there exists a smooth function $y$ defined on the whole real line 
and satisfying the differential equation $qy'+q'y=q''$. 
For instance, if $Q=x^2\frac{d}{dx}$,
the Atiyah class is non-trivial 
since no function $y\in C^\infty(\RR)$ satisfies $x^2\frac{dy}{dx}+2xy=2$ 
and therefore 
there exists no $Q$-invariant connection on $\RR$.
\end{example}

\section{Formality theorem and Kontsevich--Duflo theorem
 for $\mathfrak{g}$-manifolds}

The main results of this Note are a formality theorem for
 $\mathfrak{g}$-manifolds 
and its consequence: a Kontsevich--Duflo type theorem for $\mathfrak{g}$-manifolds. 

\begin{theorem}[Formality theorem for $\mathfrak{g}$-manifolds] 
\label{formality}
Given a $\mathfrak{g}$-manifold $M$ 
and an affine torsionfree connection $\nabla$ on $M$,
there exists an $L_\infty$ quasi-isomorphism $\Phi$ 
from the dgla
$\tot\big(\Lambda^\bullet \mathfrak{g}\dual \otimes_{\KK} \Tpoly\bullet\big)$
to the dgla
$\tot \big(\Lambda^\bullet\mathfrak{g}\dual \otimes_{\KK} \Dpoly\bullet\big)$
whose first `Taylor coefficient' $\Phi_1$ satisfies the following two properties: 
\begin{enumerate} 
\item $\Phi_1$ is, up to homotopy, 
an isomorphism of associative algebras (and hence induces an isomorphism of associative algebras of the cohomologies); 
\item $\Phi_1$ is equal to the composition 
$\hkr\circ\cTodd^{\frac{1}{2}}_{M/\mathfrak{g}}$ 
of the HKR map and the action of the square root of the Todd cocycle 
$\cTodd^{\frac{1}{2}}_{M/\mathfrak{g}} \in \bigoplus_{k=0} \Lambda^{k}\mathfrak{g}^\vee \otimes \Omega^k (M)$ on
$\tot\big(\Lambda^\bullet\mathfrak{g}\dual\otimes_{\KK}\Tpoly\bullet\big)$
by contraction. 
\end{enumerate}
\end{theorem}

As an immediate consequence, we have the following
\begin{theorem}[Kontsevich--Duflo type theorem for $\mathfrak{g}$-manifolds] 
\label{kontsevich-duflo}
Given a $\mathfrak{g}$-manifold $M$, the map
\[ \hkr\circ\Todd^{\frac{1}{2}}_{M/\mathfrak{g}}:
H^\bullet_{\CE}\big(\mathfrak{g},\Tpolycmplx\big)
\xrightarrow{\cong} H^\bullet_{\CE}\big(\mathfrak{g},\Dpolycmplx\big) \]
is an isomorphism of Gerstenhaber algebras.
It is understood that the square root $\Todd^{\frac{1}{2}}_{M/\mathfrak{g}}$ of the 
Todd class $\Todd_{M/\mathfrak{g}}\in\bigoplus_{k=0} 
H_{\CE}^{k}\big(\mathfrak{g},\Omega^k(M)\big)$ acts on 
$H^\bullet_{\CE}\big(\mathfrak{g},\Tpolycmplx\big)$ by contraction. 
\end{theorem}


Theorem~\ref{formality} follows from a more general result of ours, 
a formality theorem for Lie pairs, whose detailed proof 
will appear in a forthcoming revision of~\cite{LSX:16}. 
A pair of Lie algebroids (or Lie pair in short) consists of 
a Lie algebroid $L$ and a Lie subalgebroid $A$ of $L$. 
Given any Lie pair, our formality theorem for Lie pairs 
establishes an $L_\infty$ quasi-isomorphism $\Phi$ from 
the polyvector fields `on the pair' to the polydifferential 
operators `on the pair.'
The first `Taylor coefficient' $\Phi_1$ of the $L_\infty$ 
quasi-isomorphism $\Phi$ preserves the associative algebra 
structures up to homotopy and admits an explicit description 
in terms of the Hochschild--Kostant--Rosenberg map 
and the Todd cocycle of the Lie pair.
Now every $\mathfrak{g}$-manifold $M$ determines in a canonical way 
a matched pair: $(\mathfrak{g}\ltimes M,T_M)$ 
\cite[Example~5.5]{MR1460632} \cite{MR1650045}.
The notation $\mathfrak{g}\ltimes M$ refers to the transformation 
Lie algebroid arising from the infinitesimal 
$\mathfrak{g}$-action on $M$.
Therefore, we can form a Lie pair $(L,A)$, 
where $L=(\mathfrak{g}\ltimes M)\bowtie T_M$ 
and $A=\mathfrak{g}\ltimes M$. 
For this particular pair, the polyvector fields 
and polydifferential operators reduce to
$\tot\big(\Lambda^\bullet\mathfrak{g}\dual\otimes_{\KK}
\Tpoly\bullet\big)$ and
$\tot\big(\Lambda^\bullet\mathfrak{g}\dual\otimes_{\KK}
\Dpoly\bullet\big)$ respectively.
Theorem~\ref{formality} then follows from our formality theorem 
for Lie pairs~\cite{LSX:16}. 

To the best of our knowledge, the first construction of an $L_\infty$ quasi-isomorphism  
from the dgla $\tot\big(\Lambda^\bullet \mathfrak{g}\dual \otimes_{\KK} \Tpoly\bullet\big)$
to the dgla $\tot \big(\Lambda^\bullet\mathfrak{g}\dual \otimes_{\KK} \Dpoly\bullet\big)$
can be credited to Dolgushev~\cite[concluding remarks]{MR2102846}.

Applications of Theorem~\ref{formality} to the deformation quantization 
of $\mathfrak{g}$-manifolds will be considered elsewhere.


\section*{Acknowledgements}
The authors thank Ruggero Bandiera, Martin Bordemann, Vasily Dolgushev, 
Olivier Elchinger, Marco Manetti, Boris Shoikhet, Jim Stasheff, and Dima Tamarkin 
for inspiring discussions and useful comments.
Mathieu Sti\'enon would like to express his gratitude to the 
Universit\'e Paris--Diderot for its hospitality while work 
on this project was underway.  

\bibliographystyle{amsplain}
\bibliography{references}

\def\cprime{$'$}
\providecommand{\bysame}{\leavevmode\hbox to3em{\hrulefill}\thinspace}
\providecommand{\MR}{\relax\ifhmode\unskip\space\fi MR }
\providecommand{\MRhref}[2]{%
  \href{http://www.ams.org/mathscinet-getitem?mr=#1}{#2}
}
\providecommand{\href}[2]{#2}
\begin{thebibliography}{10}

\bibitem{MR0086359}
Michael~F. Atiyah, \emph{Complex analytic connections in fibre bundles}, Trans.
  Amer. Math. Soc. \textbf{85} (1957), 181--207. \MR{0086359}

\bibitem{MR2646112}
Damien Calaque and Michel Van~den Bergh, \emph{Hochschild cohomology and
  {A}tiyah classes}, Adv. Math. \textbf{224} (2010), no.~5, 1839--1889.
  \MR{2646112}

\bibitem{MR3439229}
Zhuo Chen, Mathieu Sti{\'e}non, and Ping Xu, \emph{From {A}tiyah classes to
  homotopy {L}eibniz algebras}, Comm. Math. Phys. \textbf{341} (2016), no.~1,
  309--349. \MR{3439229}

\bibitem{MR2102846}
Vasiliy Dolgushev, \emph{Covariant and equivariant formality theorems}, Adv.
  Math. \textbf{191} (2005), no.~1, 147--177. \MR{2102846 (2006c:53101)}

\bibitem{MR0161898}
Murray Gerstenhaber, \emph{The cohomology structure of an associative ring},
  Ann. of Math. (2) \textbf{78} (1963), 267--288. \MR{0161898}

\bibitem{MR0142598}
Gerhard Hochschild, Bertram Kostant, and Alex Rosenberg, \emph{Differential
  forms on regular affine algebras}, Trans. Amer. Math. Soc. \textbf{102}
  (1962), 383--408. \MR{0142598}

\bibitem{MR2062626}
Maxim Kontsevich, \emph{Deformation quantization of {P}oisson manifolds}, Lett.
  Math. Phys. \textbf{66} (2003), no.~3, 157--216. \MR{2062626 (2005i:53122)}

\bibitem{LSX:16}
Hsuan-{Y}i Liao, Mathieu Sti{\'e}non, and Ping Xu, \emph{Formality theorem and
  {K}ontsevich--{D}uflo type theorem for {L}ie pairs}, ArXiv e-prints (2016).

\bibitem{MR1650045}
K.~C.~H. Mackenzie, \emph{Drinfel$\prime$d doubles and {E}hresmann doubles for
  {L}ie algebroids and {L}ie bialgebroids}, Electron. Res. Announc. Amer. Math.
  Soc. \textbf{4} (1998), 74--87. \MR{1650045}

\bibitem{MR3319134}
Rajan~A. Mehta, Mathieu Sti{\'e}non, and Ping Xu, \emph{The {A}tiyah class of a
  dg-vector bundle}, C. R. Math. Acad. Sci. Paris \textbf{353} (2015), no.~4,
  357--362. \MR{3319134}

\bibitem{MR1460632}
Tahar Mokri, \emph{Matched pairs of {L}ie algebroids}, Glasgow Math. J.
  \textbf{39} (1997), no.~2, 167--181. \MR{1460632}

\bibitem{MR2699544}
Dmitry~E. Tamarkin, \emph{Operadic proof of {M}. {K}ontsevich's formality
  theorem}, ProQuest LLC, Ann Arbor, MI, 1999, Thesis (Ph.D.)--The Pennsylvania
  State University. \MR{2699544}

\end{thebibliography}

\end{document}